\documentclass[final,5p,times,twocolumn]{elsarticle}

\usepackage{amsmath}
\usepackage{amssymb}
\usepackage{xcolor}
\usepackage{todonotes}
\usepackage{ulem}
\usepackage{soul}
\usepackage{url}
\usepackage{hyperref}

\usepackage[english]{babel}
\graphicspath{{figs/}{scripts/}}

\renewcommand{\vec}[1]{\boldsymbol{#1}}
\newcommand{\jt}[1]{{#1}}

\newcounter{bla}

\journal{Computer Physics Communications}

\begin{document}

\begin{frontmatter}



\title{Geometric multigrid method for solving Poisson's equation on octree grids with irregular boundaries}


\author[a]{Jannis Teunissen\corref{author}}
\author[b]{Francesca Schiavello}

\cortext[author] {Corresponding author.\\\textit{E-mail address:} jannis.teunissen@cwi.nl}
\address[a]{Centrum Wiskunde \& Informatica, Amsterdam, the Netherlands}
\address[b]{UKRI, STFC, Hartree Centre, United Kingdom}

\begin{abstract}
  A method is presented to include irregular domain boundaries in a geometric multigrid solver.
  Dirichlet boundary conditions can be imposed on an irregular boundary defined by a level set function.
  Our implementation employs quadtree/octree grids with adaptive refinement, a cell-centered discretization and pointwise smoothing.
  Boundary locations are determined at a subgrid resolution by performing line searches.
  For grid blocks near the interface, custom operator stencils are stored that take the interface into account.
  For grid block away from boundaries, a standard second-order accurate discretization is used.
  The convergence properties, robustness and computational cost of the method are illustrated with several test cases.
\end{abstract}

\begin{keyword}
multigrid, irregular boundary, Poisson equation, adaptive mesh refinement, level set function

\end{keyword}

\end{frontmatter}



{\bf NEW VERSION PROGRAM SUMMARY}

\begin{small}
\noindent
{\em Program Title: Afivo}                                          \\
{\em CPC Library link to program files:} (to be added by Technical Editor) \\
{\em Developer's repository link:} \url{https://github.com/MD-CWI/afivo} \\
{\em Code Ocean capsule:} (to be added by Technical Editor)\\
{\em Licensing provisions(please choose one):} GPLv3 \\
{\em Programming language:} Fortran \\
{\em Journal reference of previous version:} \url{https://doi.org/10.1016/j.cpc.2018.06.018}                  \\
{\em Does the new version supersede the previous version?:} yes  \\
{\em Reasons for the new version:} Add support for internal boundaries in the geometric multigrid solver\\
{\em Summary of revisions:} The geometric multigrid solver was generalized in several ways: a coarse grid solver from the Hypre library is used, operator stencils are now stored per grid block, and methods for including boundaries via a level set function were added. \\
{\em Nature of problem:} The goal is to solve Poisson's equation in the presence of irregular boundaries that are not aligned with the computational grid. It is assumed these irregular boundaries are defined by a level set function, and that a Dirichlet type boundary condition is applied. The main applications are 2D and 3D simulations with octree-based adaptive mesh refinement, in which the mesh frequently changes but the irregular boundaries do not. \\
{\em Solution method:} A geometric multigrid method compatible with octree grids is developed, using a cell-centered discretization and point-wise smoothing. Near irregular boundaries, custom operator stencils are stored. Line searches are performed to locate interfaces with sub-grid resolution. To increase the methods robustness, this line search is modified on coarse grids if boundaries are otherwise not resolved. The multigrid solver uses OpenMP parallelization.\\


\end{small}

\section{Introduction}
\label{sec:introduction}


A common elliptic partial differential equation (PDE) is Poisson's equation
\begin{equation}
  \label{eq:poisson}
  \nabla \cdot (a(\vec{x}) \nabla \phi) = g,
\end{equation}
where the right-hand side $g$ and coefficient $a(\vec{x})$ are given and $\phi$
has to be obtained given certain boundary conditions.
Equation~(\ref{eq:poisson}) can numerically be solved with a variety of techniques,
for example using fast Fourier transforms (FFTs), cyclic
reduction, direct sparse solvers, (preconditioned) Krylov methods,
multipole methods and multigrid methods, see e.g.~\cite{Douglas_2003,Gholami_2016}.
The most suitable method depends on the type of computational grid, the boundary conditions, the spatial variation in $a(\vec{x})$, and the available computational hardware.
Our goal is to develop an efficient geometric multigrid scheme for the following case:
\begin{itemize}
  \item There are irregular Dirichlet boundary conditions.
  These boundaries are located inside the computational domain, but they are not aligned with the numerical grid.
  \item Equation~(\ref{eq:poisson}) has to be solved several times for different right-hand sides, but with the same irregular boundaries.
  \item The coefficient $a(\vec{x})$ is constant.
  \item The computational grid is a quadtree/octree mesh that is frequently adapted, so that it is desirable to have a (mostly) matrix-free method.
\end{itemize}

Multigrid methods~\cite{Hackbusch_1985,Trottenberg_2000_multigrid,Brandt_2011} can be used to solve equations like~(\ref{eq:poisson}) with great efficiency.
The main idea is to iteratively damp the error on a hierarchy of grids with a smoother.
On coarse grids, the long-wavelength components of the error are damped, and on fine grids the short-wavelength components.
Information from different grid levels is combined via prolongation (i.e., interpolation) to finer grids, and via restriction to coarser grids.
Multigrid methods can have a computational cost linear in the number of unknowns, which is ideal.
We focus on geometric multigrid (GMG) methods, which solve problems on a given hierarchy of numerical grids.
In contrast, algebraic multigrid (AMG) methods can be used to solve more general linear systems.
This flexibility is attractive for problems with irregular boundaries, but the cost of AMG methods is generally higher~\cite{Stuben_2001}.

Considerable work has been done on solving equations like~(\ref{eq:poisson}) with geometric multigrid in the presence of irregular boundaries.
We briefly mention some relevant work below.
In~\cite{Wan_2004}, a matrix-free geometric multigrid was developed that could handle irregular boundaries, which were tracked by a level-set function.
Node-centered grids were considered from 1D to 3D, and interpolation was performed by locally solving the elliptic PDE for a grid point.
Besides Dirichlet boundaries, the authors also consider discontinuities in the PDE coefficient $a(\vec{x})$.
In~\cite{Guillet_2011}, a geometric multigrid scheme was presented to apply irregular Dirichlet boundary conditions on AMR grids, with a focus on self-gravitating astrophysical flows.
On the fine grid boundaries were described by a mask, leading to a staircase pattern.
The authors discuss a trade-off between a first and second order accurate scheme for representing boundaries, with the second order scheme suffering from a lack of convergence on coarse grids when boundaries are not well resolved.
In~\cite{Theillard_2013}, a multigrid solver was presented for elliptic and parabolic problems on quadtree and octree grids.
A node-centered discretization was used and irregular boundaries were described by a level-set function.
So-called ghost values were obtained by third-order extrapolation near refinement boundaries.
In~\cite{Botto_2013}, a geometric multigrid solver compatible with irregular Neumann boundaries was presented.
The boundaries were represented by a staircase pattern on the fine grid.
The authors highlight the importance of a conservative discretization, which is also referred to as a compatibility condition, see e.g.~\cite{Lee_2007,Teunissen_2018}.
Interpolation was avoided near boundaries, leading to a first order accurate method.
In~\cite{Weber_2015}, a cut-cell geometric multigrid solver was presented supporting both Dirichlet and Neumann boundary conditions, with a focus on the efficient simulation and visualization of incompressible flow.
A cell-centered discretization was used, and a constant (zeroth-order) prolongation scheme.
The method was shown to be first order accurate for Dirichlet boundaries and second order accurate for Neumann boundaries.

The main novelty of the method presented here is that it combines the following aspects:
\begin{itemize}
  \item The flexible handling of different geometries via a level-set function.
  \item An (approximately) second-order accurate cell-centered discretization that is compatible with adaptive mesh refinement (AMR) on quadtree/octree grids.
  \item The use of a line search method to accurately locate interfaces.
  \item A correction for unresolved boundaries on coarse grids.
  \item An efficient open-source implementation, with custom stencils only stored for grid blocks that contain a boundary.
\end{itemize}

\section{Multigrid method without irregular boundaries}
\label{sec:method-no-boundaries}

Below, the basis of multigrid method used in this paper is briefly introduced.
The extension to irregular boundaries is discussed in section~\ref{sec:method-boundaries}.

\subsection{Mesh}
\label{sec:mesh}

\begin{figure}
  \centering
  \includegraphics[width=8cm]{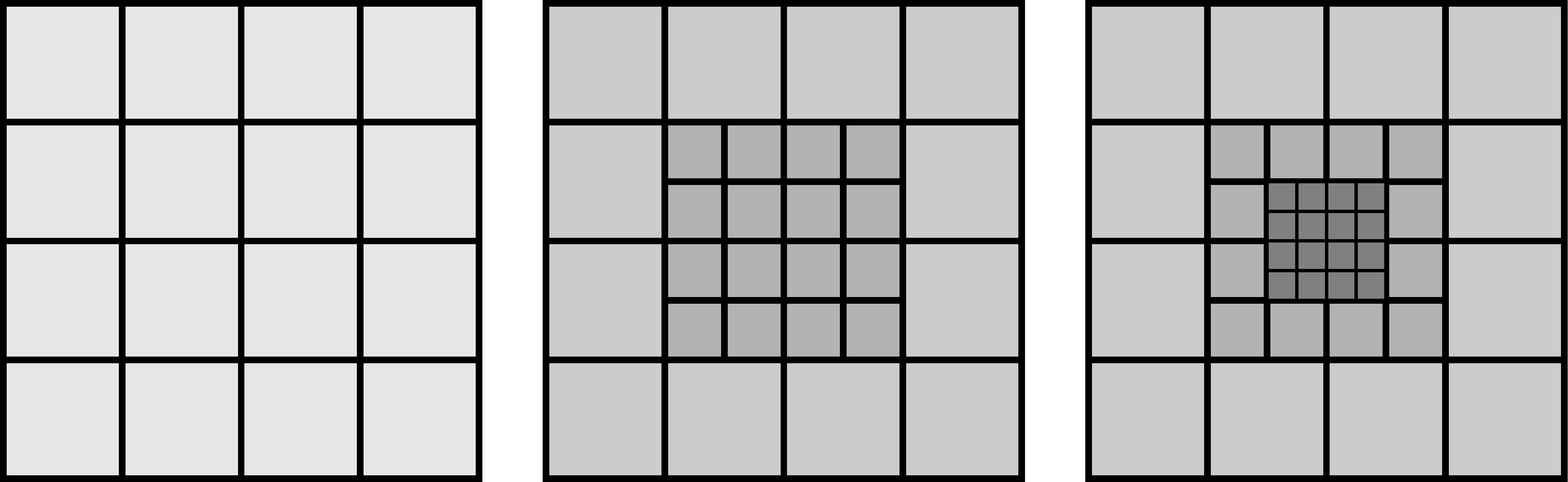} %
  \caption{Illustration of a quadtree grid.
    The squares indicate grid blocks, which each contain $N \times N$ cells (not indicated).
    From left to right, refinement is added around the center.}
  \label{fig:quadtree}
\end{figure}

We consider so-called octree meshes, see figure~\ref{fig:quadtree}.
In our implementation, which is based on the \texttt{afivo} framework~\cite{Teunissen_2018}, such a mesh consists of blocks of $N^D$ cells, where $D$ denotes the problem dimension.
These blocks can be refined by halving the grid spacing, so that $2^D$ refined child blocks cover a parent block.
Nearby blocks are refined, if necessary, to ensure that adjacent blocks differ by at most one refinement level.
A tree fulfilling such a condition is called 2:1 balanced.

Octree meshes balance adaptivity and computational efficiency.
Because each block has the same shape, computations, communication and mesh refinement can be implemented rather efficiently.
We use a cell-centered approach, in which the solution $\phi$ and right-hand side of
equation~(\ref{eq:poisson}) are defined at cell centers, and the components of $\nabla \phi$ are defined at cell faces.

\subsection{Geometric multigrid method}
\label{sec:multigrid}

The algorithms presented here for irregular boundaries extend the geometric multigrid solver of the Afivo framework~\cite{Teunissen_2018}.
This solver implements the Full Approximation Scheme (FAS), in which the solution $\phi$ is approximated on all grid levels.
A brief overview is given below; further details can be found in~\cite{Teunissen_2018} and in~\cite{Teunissen_2019}, which describes an MPI-parallel version.

\textbf{Operators} A standard finite difference discretizations of the Laplacian is used, with 3, 5 and 7-point numerical stencils in 1D, 2D and 3D, respectively.
If the grid spacing $\Delta x$ is constant and there are no boundaries, a second order accurate discretization of equation~(\ref{eq:poisson}) in 1D is given by
\begin{equation}
  \label{eq:lpl-1d}
  \nabla^2 \phi_i = \frac{1}{\Delta x} \left(\frac{\phi_{i+1} - \phi_{i}}{\Delta x}
    - \frac{\phi_{i} - \phi_{i-1}}{\Delta x}\right) = g_i,
\end{equation}
where $g$ is the right-hand side.
The residual for an approximate solution $\phi^*$ is defined as
\begin{equation}
  \label{eq:residual}
  r = g - \nabla^2 \phi^*.
\end{equation}

\textbf{Smoother} Gauss-Seidel red-black (GSRB) smoothers are used.
The unknowns are first divided into red and black groups, in a checkerboard fashion.
Equations like~(\ref{eq:lpl-1d}) can then be solved in parallel for one group, assuming the other group's values stay fixed.
For equation~(\ref{eq:lpl-1d}), this results in
\begin{equation*}
  \label{eq:gsrb-1d}
  \phi_{i} = \frac{1}{2} \left(\phi_{i+1} + \phi_{i-1} + \Delta x^2 g\right).
\end{equation*}
The error is damped by alternatingly solving for the red and black groups.

\textbf{Prolongation and restriction} Prolongation is the transfer of coarse-grid corrections to a finer grid, which is a key part of a geometric multigrid method.
Standard (bi/tri)linear interpolation is here used, also when irregular boundaries are present.
Restriction is the transfer of information to a coarser grid.
This is implemented by taking the average of the $2^D$ fine-grid cells covering a coarse grid cell.

\begin{figure}
  \centering
  \includegraphics[width=8cm]{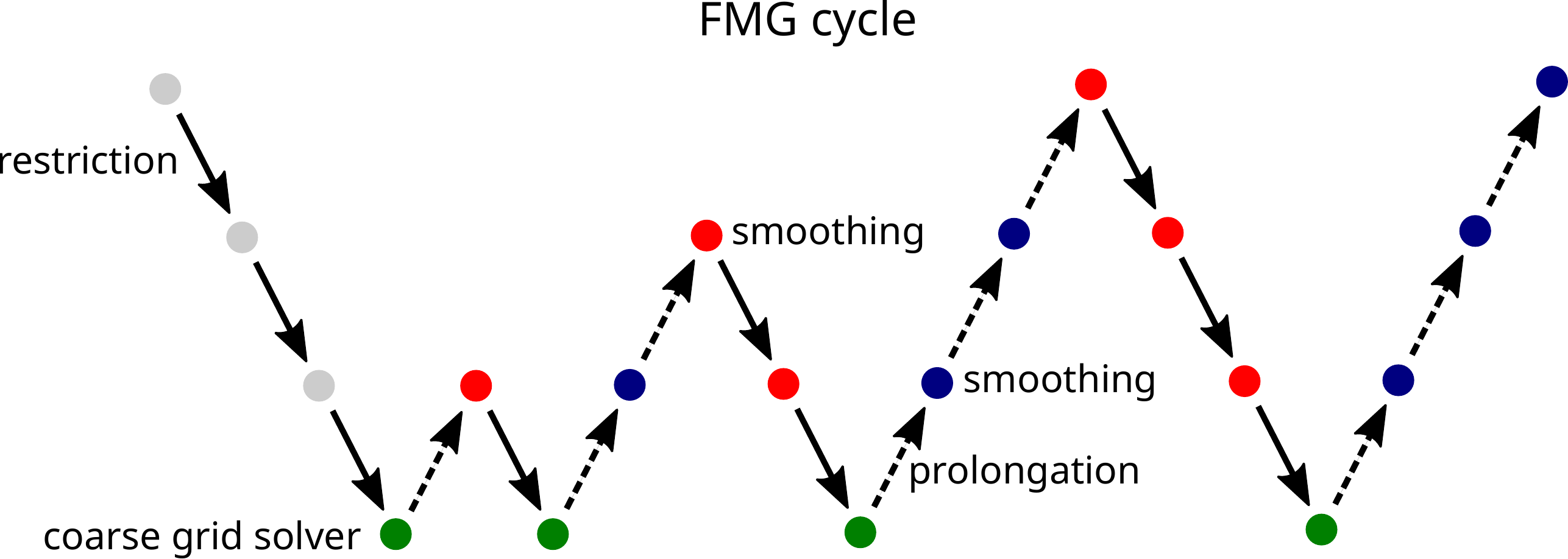}
  \caption{Illustration of the FMG cycle for a grid with four levels.}
  \label{fig:mg-cycles}
\end{figure}

\textbf{Multigrid cycle} A standard V-cycle and full multigrid (FMG) cycle are implemented, see figure~\ref{fig:mg-cycles}.
The V-cycle goes from fine to coarse, and then back to fine.
The FMG cycle iteratively performs V-cycles from the coarsest grid up to the finest grid.
Although FMG cycles are more expensive than V-cycles, they can guarantee a reduction of the residual that is independent of the problem size~\cite{Trottenberg_2000_multigrid}.
This results in the ideal $\mathcal{O}(N)$ computational cost of FMG, where $N$ is the number of unknowns.

Every time a grid level is visited in a cycle, smoothing is performed.
In the upward part of a cycle, $N_\mathrm{up}$ smoothing steps are performed before prolongation, and in the downward part of a cycle, $N_\mathrm{down}$ smoothing steps performed before restriction.
We here use $N_\mathrm{up} = N_\mathrm{down} = 2$.
The handling of the coarse grid is discussed in section~\ref{sec:coarse-grid-solver}.

\textbf{Ghost cells} When performing multigrid on an adaptive mesh, it is convenient to extend grid blocks with a layer of ghost cells.
It is important that the ghost cells near refinement boundaries are filled in such a way that the fine-grid discretization is consistent with the underlying coarse grid.
We here follow the same approach as in~\cite{Teunissen_2018}.
The basic idea is that ghost cells are filled in such a way that the coarse and averaged fine `flux' across the refinement boundary (e.g., $\partial_x \phi$) agree.

\section{Implementation of boundaries}
\label{sec:method-boundaries}

\subsection{Level set function}
\label{sec:level-set-function}

Internal boundaries are here defined by the zero contour of a level set function (LSF)~\cite{Osher_1988,Sethian_1999}:
\begin{equation*}
  \label{eq:lsf}
  f(\vec{x}) = 0.
\end{equation*}
For example, a spherical boundary of radius $R$ centered at $\vec{x}_c$ can be
described by
\begin{equation}
  \label{eq:lsf-example}
  f(\vec{x}) = ||\vec{x} - \vec{x}_c|| - R.
\end{equation}
In this case, the LSF is a signed distance function, with a negative sign inside the sphere.
Examples of other LSFs are given in section~\ref{sec:num-experiments}.



\subsection{Distance computation}
\label{sec:lsf-distance}

For geometric multigrid, it is important that the locations of boundaries (i.e., roots of the LSF) agree well between grid levels.
We therefore use a line search method to locate boundaries at a sub-grid resolution.

Let $\vec{a}$ denote a start point, e.g., the center of a grid cell, and $\vec{b}$ a neighboring point.
We want to know if there is a root in the LSF on the line segment from $\vec{a}$ to $\vec{b}$, and if so, how far this root is from $\vec{a}$.
This information is here stored in a single value $d$, which denotes the relative distance to the boundary.
If there is no boundary between $\vec{a}$ and $\vec{b}$, $d = 1$.
Otherwise, if there is root at $\vec{x}_0$, $d$ is given by
\begin{equation}
  \label{eq:relative-distance}
  d = ||\vec{x}_0 - \vec{a}|| / ||\vec{b} - \vec{a}||.
\end{equation}

\begin{figure*}
  \centering
  \includegraphics[width=14cm]{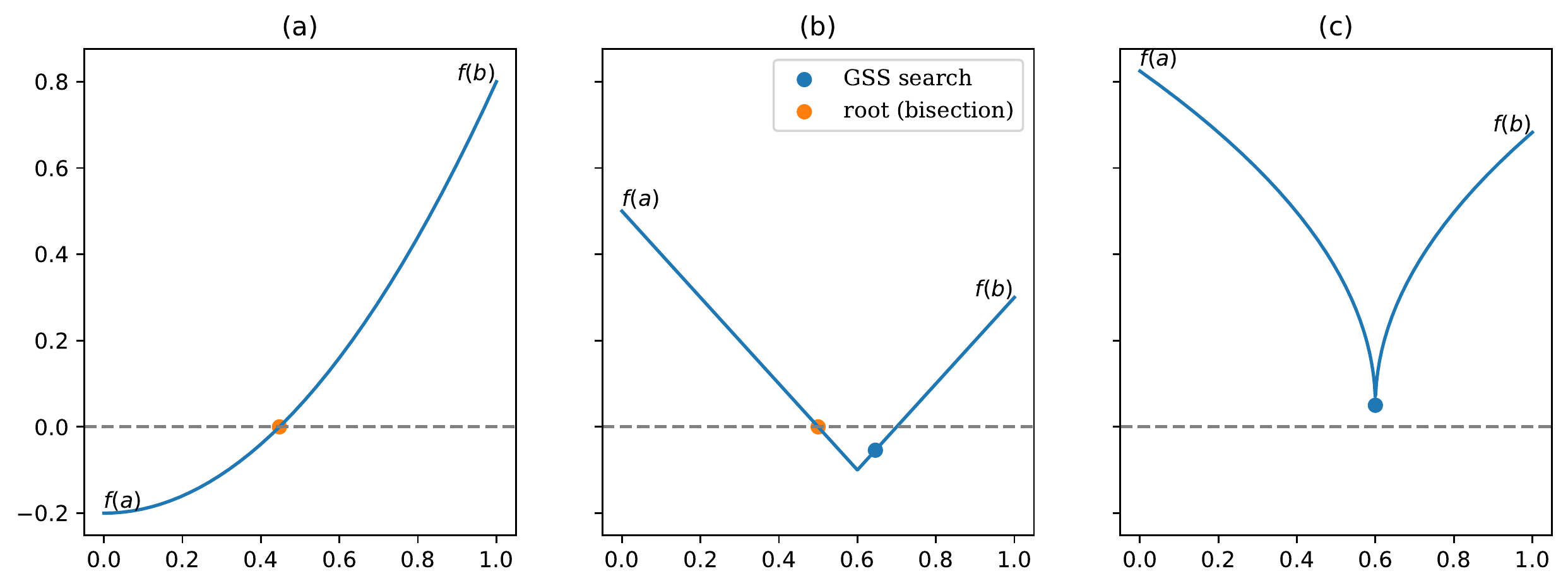}
  \caption{Illustration of the procedure for locating roots in the level set function between two points $\vec{a}$ and $\vec{b}$. Three cases are illustrated: (a) The root is already bracketed, so bisection is directly applied; (b) A bracket is determined using Golden section search, then bisection is applied; (c) No bracket is found, and thus also no root.}
  \label{fig:rootfinding}
\end{figure*}

The procedure for locating roots is illustrated in figure~\ref{fig:rootfinding}.
If $f(\vec{a}) \times f(\vec{b}) \leq 0$, bisection is used to locate the root $\vec{x}_0$ between $\vec{a}$ and $\vec{b}$, with a relative tolerance of $\epsilon_\mathrm{tol}$.
Otherwise, a bracket for the potential root first has to be determined.
We use golden section search to minimize $f(\vec{x}) \times f(\vec{a})$ on the line from $\vec{a}$ to $\vec{b}$.
As soon as $f(\vec{x}) \times f(\vec{a}) \leq 0$, bisection is again applied on the interval between $\vec{a}$ and $\vec{x}$.
If this condition is not met within a given number of iterations, corresponding to the same relative tolerance $\epsilon_\mathrm{tol}$, it is assumed there is no boundary.

Note that if there are two roots on the interval between $\vec{a}$ and $\vec{b}$, the bracket search will eliminate the one farthest from $\vec{a}$.
If there are three or more roots, it is not guaranteed that the above procedure finds root closest to $\vec{a}$.
By default, we use a small relative tolerance of $\epsilon_\mathrm{tol} = 10^{-8}$.

\subsection{Discretization of Laplacian with boundaries}
\label{sec:laplacian-operator}

When there are irregular boundaries, the distances between an unknown $\phi_{i}$
and neighboring values are no longer fixed. The numerical Laplacian of
equation~(\ref{eq:lpl-1d}) can then be generalized to
\begin{equation}
  \label{eq:lpl-1d-lsf}
  \nabla^2 \phi_i = \frac{2}{(d_{i+1} + d_{i-1})\Delta x}
  \left(\frac{\phi_{i+1} - \phi_{i}}{d_{i+1} \Delta x}
    - \frac{\phi_{i} - \phi_{i-1}}{d_{i-1} \Delta x}\right) = g_i,
\end{equation}
where $0 < d_{j} \leq 1$ denotes the relative distances from $\phi_{i}$ to neighboring values $\phi_{j}$, see section~\ref{sec:lsf-distance}.
Note that in the above notation the $\phi_{j}$ (for $j \neq i$) do not always correspond to unknowns on the grid.
For example, if there is a boundary between cell $i$ and $i+1$, then $\phi_{i+1}$ will correspond to a boundary value $\phi_b$.
If the corresponding term is moved to the right-hand side, equation~(\ref{eq:lpl-1d-lsf}) becomes
\begin{equation*}
  \label{eq:lpl-1d-lsf-rhs-ex}
  \frac{2}{(d_{i+1} + d_{i-1})\Delta x}
  \left(\frac{-\phi_{i}}{d_{i+1} \Delta x}
    - \frac{\phi_{i} - \phi_{i-1}}{d_{i-1} \Delta x}\right) =
  g_i - \frac{2\phi_{b}}{(d_{i+1} + d_{i-1})d_{i+1} \Delta x^2}.
\end{equation*}

Equation~(\ref{eq:lpl-1d-lsf}) is a non-symmetric discretization that was used before in e.g.~\cite{Chen_1997,Gibou_2002,Udaykumar_1999}.
Near boundary points this discretization is first order ($\mathcal{O}(\Delta x)$) accurate, because the second derivative is not evaluated at the center of the two gradient terms.
However, if the number of boundary cells is small the global error can still be approximately second order accurate~\cite{Gibou_2002,Udaykumar_1999}.

The extension of equation~(\ref{eq:lpl-1d-lsf}) to multiple
dimensions is straightforward, with the same type of terms appearing for each
dimension. For example, in 2D, the Laplacian can be written as
\begin{eqnarray}
  \label{eq:lpl-2d-lsf}
  \nabla^2 \phi_{i,j} &= \frac{2}{(d_{i+1, j} + d_{i-1, j})\Delta x}
                        \left(\frac{\phi_{i+1, j} - \phi_{i, j}}{d_{i+1, j} \Delta x}
                        - \frac{\phi_{i, j} - \phi_{i-1, j}}{d_{i-1, j} \Delta x}\right) +
                        \nonumber \\
                      & \frac{2}{(d_{i, j+1} + d_{i, j-1})\Delta x}
                        \left(\frac{\phi_{i, j+1} - \phi_{i, j}}{d_{i, j+1} \Delta x}
                        - \frac{\phi_{i, j} - \phi_{i, j-1}}{d_{i, j-1} \Delta x}\right) = g_{i,j}.
\end{eqnarray}



\subsection{Prolongation}
\label{sec:prolongation-lsf}

Standard (bi/tri)linear prolongation is used, also when boundaries are present.
We did experiment with a custom prolongation scheme, in which a linear function was constructed between the nearest $D+1$ neighbors and/or boundaries, but this scheme did typically not lead to faster convergence.

\subsection{Coarse grid solver}
\label{sec:coarse-grid-solver}

At the coarsest grid there are essentially two options.
The first is to apply the same smoother as on other grid levels.
However, depending on the size and geometry of the coarse grid, it could take a large number of smoothing steps to achieve a desired reduction of the residual.
Therefore, we here solve the coarse grid equations using a different multigrid solver, provided by the Hypre library~\cite{Falgout_2002}.

The coarse grid is frequently visited in an FMG cycle, see figure~\ref{fig:mg-cycles}.
It is therefore important to keep the computational cost of the coarse grid solver as low as possible.
For this reason, we by default use Hypre's PFMG solver, which is a parallel semicoarsening multigrid solver that uses pointwise smoothing~\cite{Ashby_1996,Falgout_2000}.
In 1D, the PFMG solver is not available and we use Hypre's PCG solver instead.
Hypre's default tolerance of $10^{-6}$ is used for these solvers.

\subsection{Implementation aspects}
\label{sec:impl-aspects}

Below, we provide information on implementation aspects relevant for the
computational efficiency of the method.

\subsubsection{Boundary detection}
\label{sec:boundary-detection}

The line search for boundaries described in section~\ref{sec:lsf-distance} can be expensive.
For computational efficiency, we only perform such a search for cells that are sufficiently close to a boundary.
The distance to a boundary can be approximated by the ratio $|f(\vec{x})|/|| \nabla f(\vec{x}) ||$, which for a linear function $f(\vec{x})$ would be exact.
This inspires the following condition for a potential boundary, evaluated at every cell center:
\begin{equation}
  \label{eq:boundary-mask}
  |f(\vec{x})| < L \times || \nabla f(\vec{x}) ||,
\end{equation}
where $L$ should be proportional to the grid spacing $\Delta x$. For the tests
presented in section~\ref{sec:num-experiments} we use
$L = 1.5 \times \sqrt{D} \Delta x$, with $D$ the problem dimension. The
components of $\nabla f(\vec{x})$ are computed numerically using central
differencing.

Note that if $||\nabla f(\vec{x})||$ varies rapidly near boundaries, a larger safety factor than the $1.5 \times \sqrt{D}$ used above might be necessary.
One way to avoid this is to use a LSF that is (approximately) a signed distance function, such as equation~(\ref{eq:lsf-example}).

\subsubsection{Storing stencils and boundary information}
\label{sec:storing-stencils}

In our implementation, boundary information and numerical stencils are stored
per grid block. For grid blocks without boundaries, stencils are constant, so
they can be stored compactly. If there is a boundary passing through the grid
block, the following information is stored per grid cell:
\begin{itemize}
  \item The relative distances $d_i$ to boundaries from the cell center to
  neighboring cell centers, according to equation~(\ref{eq:relative-distance}).
  These distances are only stored for cells with at least one adjacent boundary.
  \item Operator stencil coefficients. For example, in 2D, five values have to
  be stored per cell for the operator in equation~(\ref{eq:lpl-2d-lsf}).
  Furthermore, the sum of the stencil coefficients that were moved to the
  right-hand side is stored, so that the value imposed at the boundary can be changed.
\end{itemize}

\subsubsection{Unresolved LSF roots on coarse grids}
\label{sec:lsf-coarse-grids}

\begin{figure}
  \centering
  \includegraphics[width=5cm]{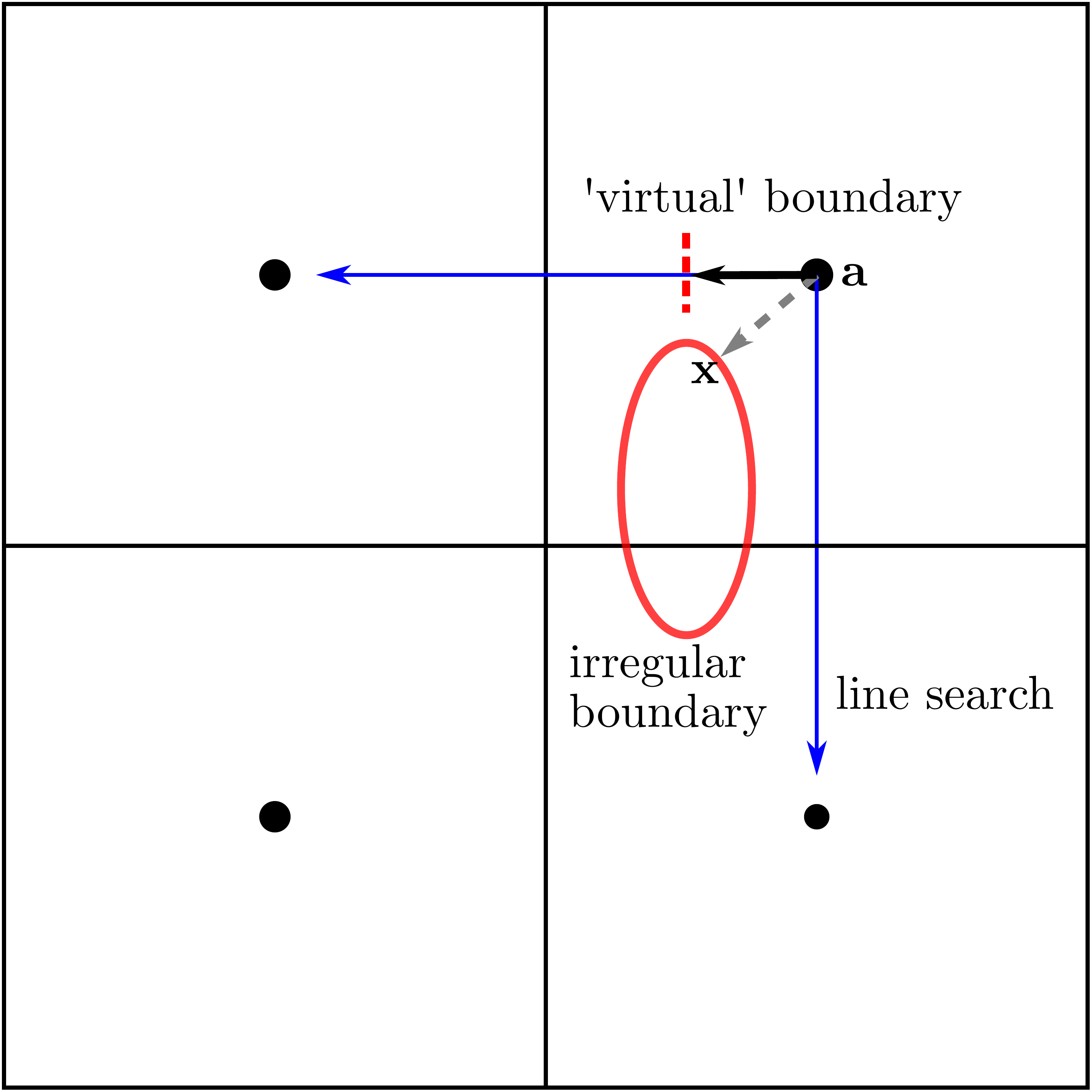}
  \caption{Illustration of procedure for unresolved boundaries on coarse grids, for the grid cell marked \textbf{a}.
    The object (red ellipse) has a small extent in both dimensions, and is therefore not detected by a line search between cell centers (blue arrows).
    First, the distance to the closest point \textbf{x} on the boundary is determined (gray arrow), and then a virtual boundary is placed between \textbf{a} and the neighboring cell center closest to \textbf{x}.
  }
  \label{fig:unresolved-boundary}
\end{figure}

If an irregular boundary has a small spatial extent in two of its dimensions, it might not be detected on a coarse grid by the line search described in section~\ref{sec:boundary-detection}.
An example is shown in figure~\ref{fig:unresolved-boundary}.
A boundary that is not detected on the coarse grid can lead to convergence issues.
We therefore perform an additional boundary search for grid cells that satisfy the following conditions:
\begin{itemize}
  \item Equation~(\ref{eq:boundary-mask}) holds, indicating there is a nearby boundary.
  \item No boundaries are detected between the cell's center and the neighboring cell centers.
  \item The grid spacing $\Delta x$ is larger than a user-defined threshold $w_\mathrm{min}$.
\end{itemize}
For these cells, gradient descent is performed in the direction in which $f$ goes to zero, starting from the cell's center ($\vec{a}$).
At most $\Delta x/w_\mathrm{min}$ steps are performed with a step size $w_\mathrm{min}$.
If after one of these steps, a location $\vec{x}$ is found such that $f(\vec{a}) f(\vec{x}) \leq 0$, a line search is performed between $\vec{a}$ and $\vec{x}$ to determine the relative distance $d$ to the boundary (normalized to the grid spacing).
This relative distance is then used in the direction of the neighboring cell closest to $\vec{x}$, as illustrated in figure~\ref{fig:unresolved-boundary}.
The resulting discretization does not accurately represent the unresolved object, but this discretization is only used on coarse grids, so it does not affect the converged fine-grid solution.


\section{Numerical experiments}
\label{sec:num-experiments}

\subsection{Convergence tests on sphere}
\label{sec:convergence-sphere}

To test the numerical convergence of the method, we solve the Laplace equation
\begin{equation*}
  \label{eq:laplace}
  \nabla^2 \phi = 0
\end{equation*}
for a spherical LSF of the form
\begin{equation}
  \label{eq:convergence-lsf}
  f_\mathrm{sphere}(\vec{x}) = ||\vec{x}|| - R,
\end{equation}
using a computational domain of unit size (e.g., the unit cube in 3D) centered at the origin.
On the spherical boundary, a Dirichlet condition $\phi = \phi_b$ is imposed.
On the boundaries of the computational domain, the following analytic solutions are imposed
\begin{align}
  \label{eq:solution-sphere}
  &\phi_\mathrm{2d}(r) = \phi_b + a \log(||\vec{x}||/R),\\
  &\phi_\mathrm{3d}(r) = \phi_b + a (1 - R/||\vec{x}||),
\end{align}
using $\phi_b = 0$ and $a = 1$.

\begin{figure}
  \centering
  \includegraphics[width=8cm]{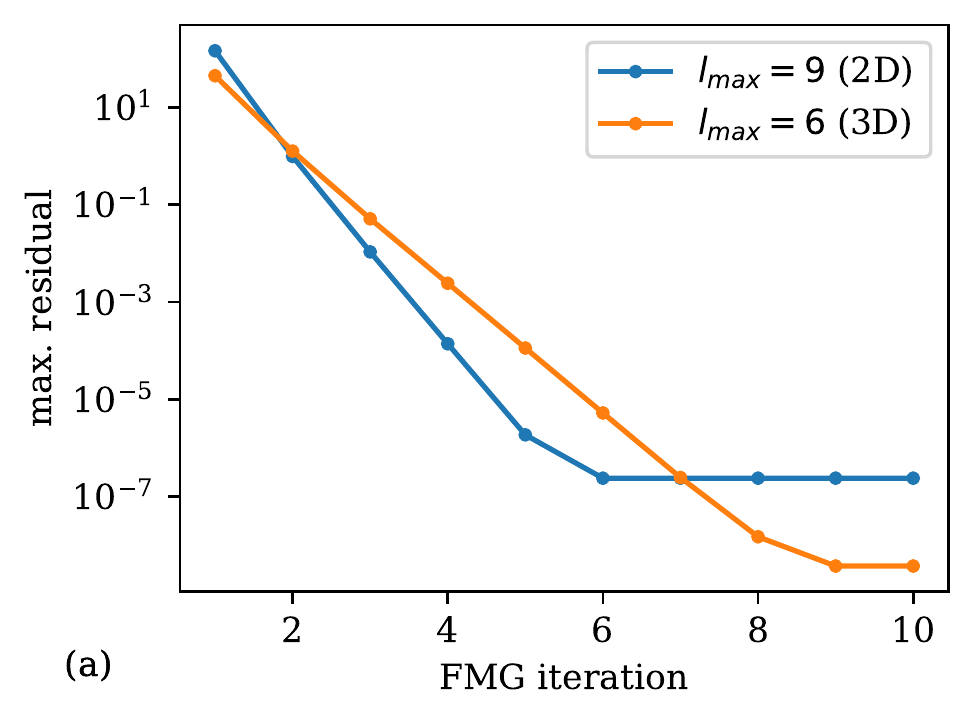}
  \includegraphics[width=8cm]{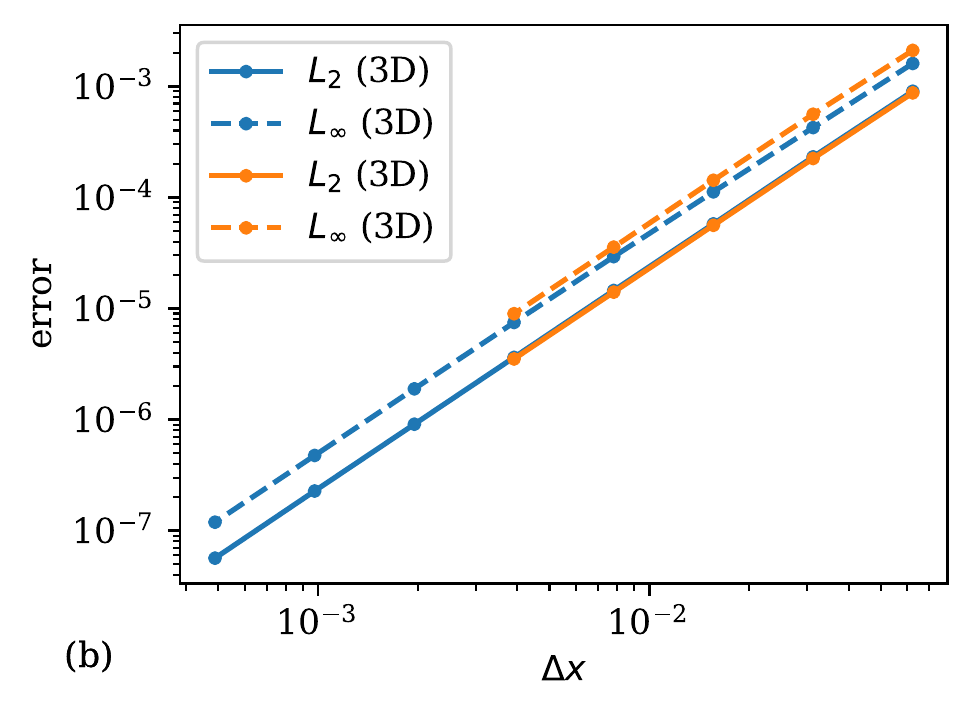}
  \caption{Convergence results the spherical LSF given by equation \eqref{eq:convergence-lsf}, with $R = 1/4$ and a uniform grid.
    a) Maximal residual versus FMG iteration, with the maximum refinement level indicated by $l_\mathrm{max}$.
    b) Error in converged solution for different grid spacings $\Delta x$, determined by comparing with equation \eqref{eq:solution-sphere}.
    Both the maximal error (indicated by $L_\infty$) and the RMSE (indicated by $L_2$) are shown.
  }
  \label{fig:sphere-convergence}
\end{figure}

As a first test, we consider the case $R = 1/4$ on uniformly refined grids.
The coarsest grid contains $8^D$ cells, and the finest grid $(8 \times 2^{l_\mathrm{max}-1})^D$ cells, where $l_\mathrm{max}$ is the maximal refinement level.
Figure~\ref{fig:sphere-convergence}a shows convergence results in 2D and 3D.
The residual reduction factor per FMG iteration is about 40--80 in 2D and about 30--40 in 3D.
Due to numerical round-off errors, the residual eventually stops decreasing.
The resulting `converged' residual is larger on finer grids because of the division by $\Delta x^2$ in equation~(\ref{eq:residual}).

After two FMG iterations, the solution error (as compared with the analytic solutions) hardly changes anymore.
For example, for the 3D case with $l_\mathrm{max} = 6$, the maximal error after the first three iterations is $0.32\times 10^{-3}$, $0.11\times 10^{-3}$ and $0.11\times 10^{-3}$.
This means that after two iterations, the discretization error dominates the convergence error.
Figure~\ref{fig:sphere-convergence}b shows that the discretization error reduces proportional to $\Delta x^2$, both in the $L_\infty$ and in the $L_2$ norm, indicating second order convergence.

\begin{figure}
  \centering
  \includegraphics[width=8cm]{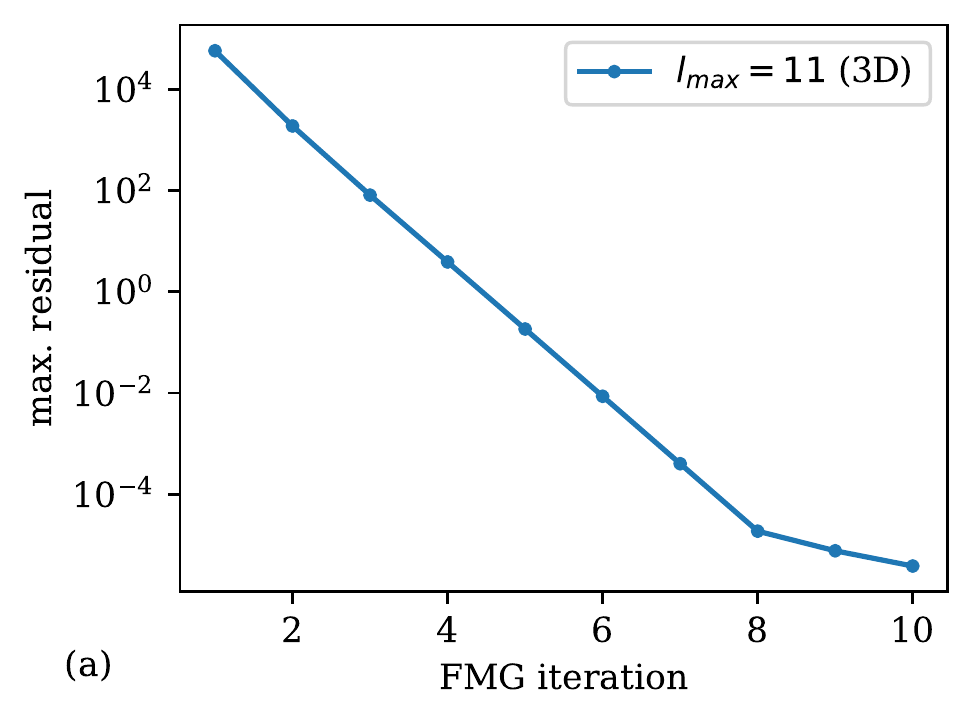}
  \includegraphics[width=8cm]{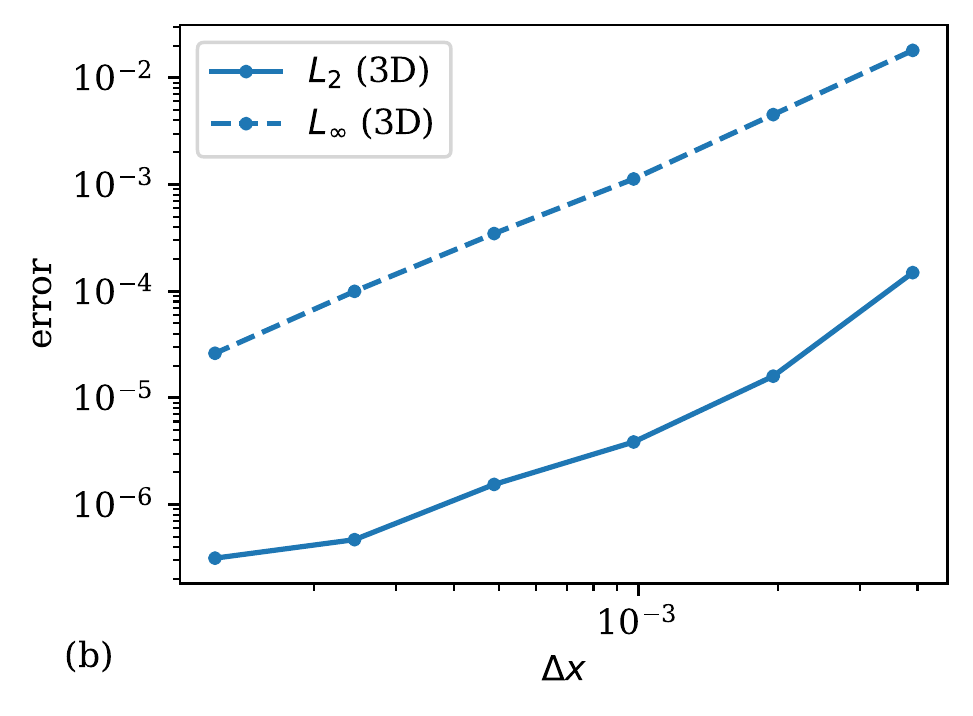}
  \caption{Convergence results for 3D spherical problem with grid refinement, in which a boundary condition is imposed at a small radius ($R = 5\times 10^{-3}$).
    a) Residual reduction. b) Error in the converged solution for different finest-grid spacings $\Delta x_\mathrm{min}$.
  }
  \label{fig:sphere-convergence-ref}
\end{figure}

As a second test, we consider the case $R = 5 \times 10^{-3}$ in 3D, in combination with grid refinement.
The following refinement criterion is used: refine if $\Delta x > \Delta x_\mathrm{min}\times\max(1, r/R)$, where $r = \sqrt{x^2 + y^2 + z^2}$ and $\Delta x_\mathrm{min}$ is the grid spacing at level $l_\mathrm{max}$.
Due to its small radius, the spherical boundary will not be resolved on the coarsest grids, but the approach described in section~\ref{sec:lsf-coarse-grids} ensures that the method still converges.
Figure~\ref{fig:sphere-convergence-ref} shows that the residual reduction factor per FMG iteration is again about 30--40.
The error in the solution is still approximately proportional to $\Delta x^2$, where $\Delta x$ is the finest grid spacing.
Note that this convergence behavior also depends on how well the mesh refinement is adapted to the problem.


\subsection{Sharp boundaries}
\label{sec:sharp-boundaries}

\begin{figure*}
  \centering
  \includegraphics[width=16cm]{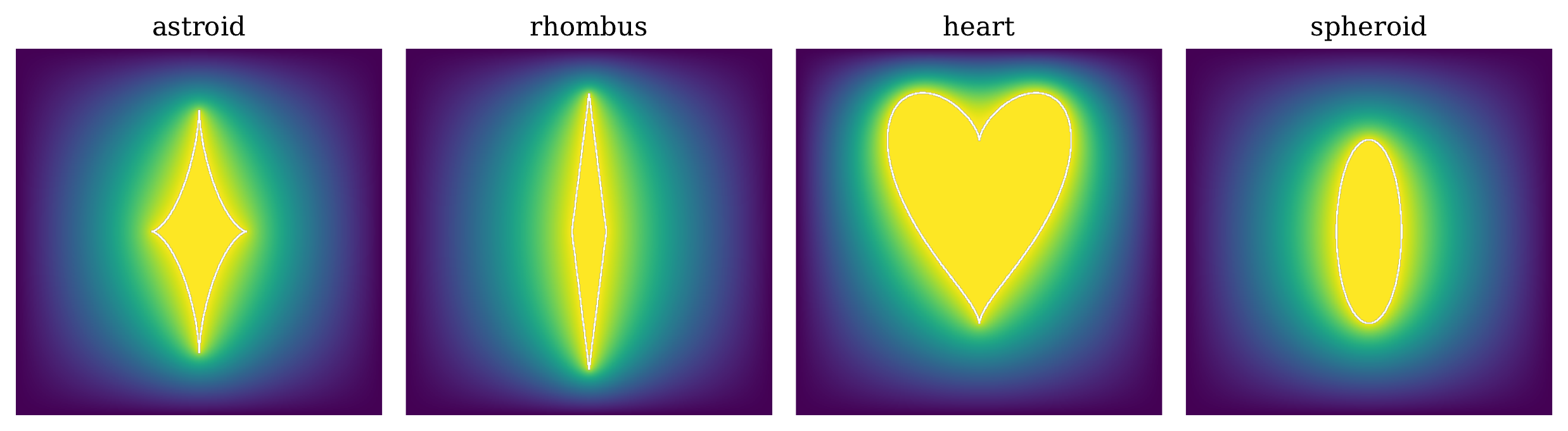}
  \caption{Illustration of solutions with the level set functions of equations~(\ref{eq:spheroid})--~(\ref{eq:astroid}), shown on the unit square.
    Contours of the boundaries are indicated in white.
    Note that the astroid shape has very sharp features.}
  \label{fig:conv-shapes}
\end{figure*}

Irregular boundaries with sharp features are a more challenging test for a multigrid-based solver.
On coarser grids, such sharp features cannot be accurately described, potentially reducing the effectiveness of the coarse-grid correction.
Furthermore, near sharp features the solution will have steep gradients, which increases interpolation errors.
To test the robustness of our solver, we consider the following level set functions
\begin{align}
  f_\mathrm{spheroid}(p, q) &= \sqrt{8 p^2 + q^2} - 1, \label{eq:spheroid}\\
  f_\mathrm{rhombus}(p, q) &= 8 |p| + |q| - 1.5,\\
  f_\mathrm{heart}(p, q) &= p^2 + (q - |p|^{2/3})^2 - 1,\\
  f_\mathrm{astroid}(p, q) &= |p|^{2/3}/0.8 + |q|^{2/3}/1.5 - 0.8. \label{eq:astroid}
\end{align}
These LSFs are evaluated on the unit square using transformed coordinates $p = (x - 0.5)/4$ and $q = (y - 0.5)/4$, see figure~\ref{fig:conv-shapes}.

As a first test, we consider uniformly refined grids in 2D of size $1024^2$ and $2048^2$, using a coarse grid size of $8\times 8$.
The corresponding maximum refinement levels are thus $l_\mathrm{max} = 8$ and $l_\mathrm{max} = 9$.
At the irregular boundary, a boundary condition $\phi = 1$ is applied, and $\phi = 0$ on the boundaries of the computational domain.
The smallest width to resolve on coarse grids (see section~\ref{sec:lsf-coarse-grids}) was set to $w_\mathrm{min} = 4\times 10^{-3}$.

\begin{figure}
  \centering
  \includegraphics[width=8cm]{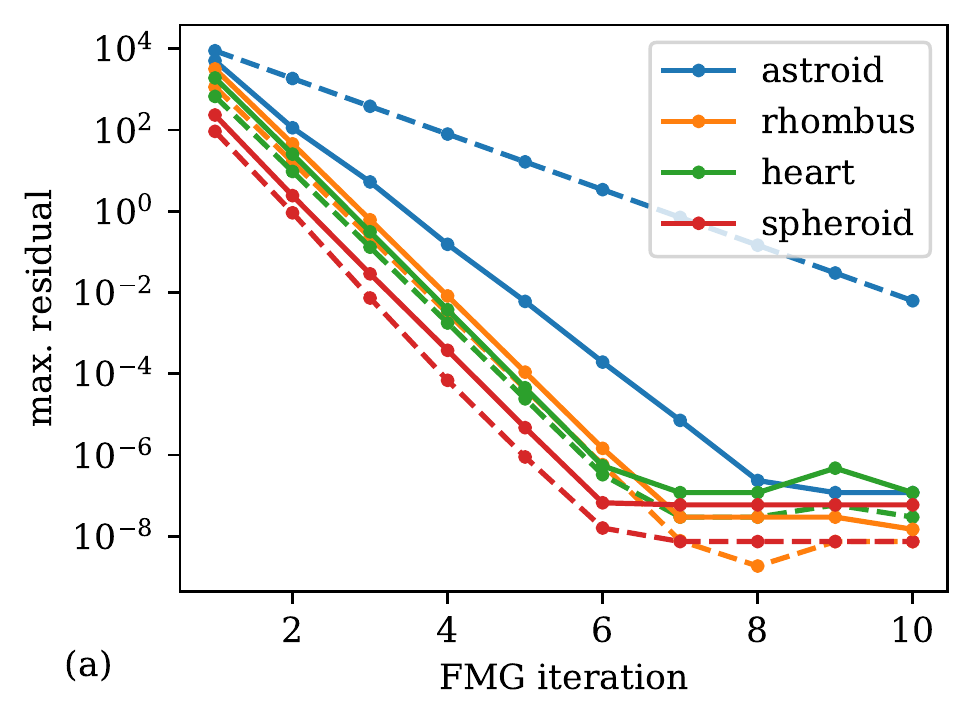}
  \includegraphics[width=8cm]{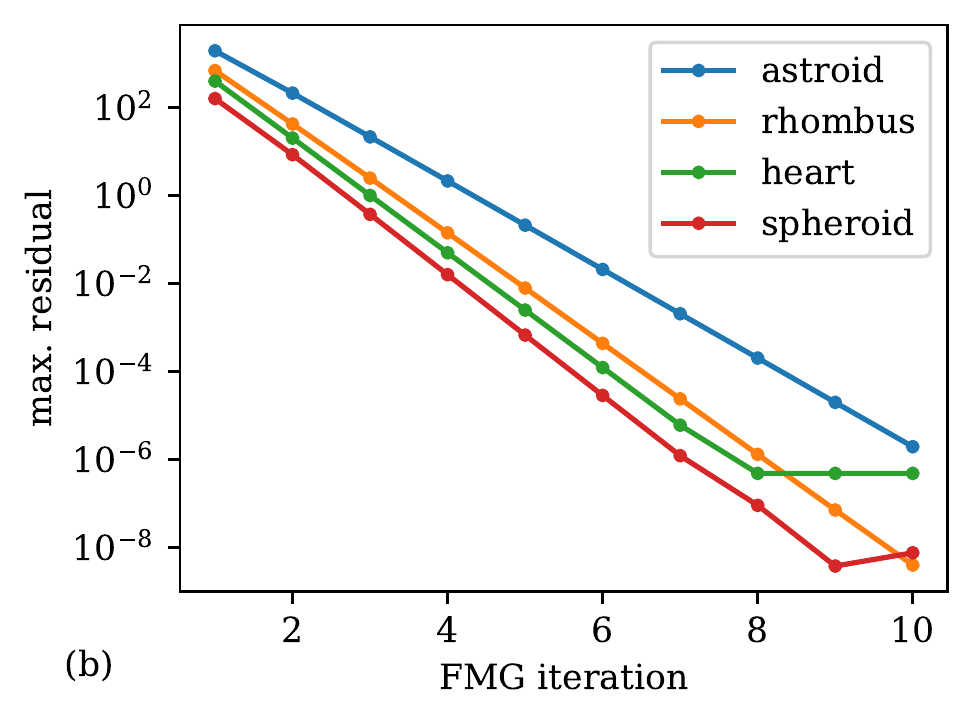}
  \caption{
    a) Maximal residual versus FMG iteration for the 2D shapes shown in figure~\ref{fig:conv-shapes}.
    The dashed lines correspond to a $1024^2$ grid, and the solid lines to a $2048^2$ grid on the unit square.
    b) Results for the axisymmetric 3D generalizations of the shapes on a $256^3$ grid.}
  \label{fig:conv-results}
\end{figure}

Figure~\ref{fig:conv-results}a shows the reduction in the residual per FMG iteration for each test case.
The residual reduction factor is similar for the spheroid, rhombus and heart shapes.
For the astroid shape the reduction factor is lower, and it changes with the refinement level.
We have also noticed that the reduction factor for this test case can depend on the position of the astroid.
This is probably due to an inconsistent description of the sharp endpoints on different refinement levels.

We generalize the above LSFs to cylindrical geometries in 3D by using transformed coordinates $p = (\sqrt{x^2 + y^2} - 0.5)/4$ and $q = (z - 0.5)/4$.
Note that in 3D there is curvature along an extra coordinate, so that solution gradients become even steeper near sharp features.
We consider the shapes given above in a 3D unit cube, using the same boundary conditions as in 2D, a coarse grid size of $8^3$ and $w_\mathrm{min} = 8\times 10^{-3}$.
Figure~\ref{fig:conv-results}b shows the residual reduction per FMG iteration on a uniformly refined grid of $256^3$ ($l_\mathrm{max} = 6$).
Note that the residual reduction factor is again lowest for the astroid shape, which has the sharpest features.

\subsection{Application example}
\label{sec:two-electrodes}

We briefly present an example relevant for the simulation of pulsed electric discharges such as streamers~\cite{Nijdam_2020}.
In discharge simulations, electrostatic fields have to be computed at every time step.
This is done by solving equation \eqref{eq:poisson} for a given electrode configuration, after which the electric field is obtained as $\vec{E} = -\nabla \phi$.
Solving Poisson's equation is typically one of the most expensive components of these simulations.
Various methods have been used to incorporate electrodes, ranging from a simple charge simulation technique (see e.g.~\cite{Luque_2008a}) to the ghost fluid method~\cite{Celestin_2009a}.
Electrodes have also been included on structured grids with different AMR framework~\cite{Kolobov_2012,Marskar_2019a}, in combination with multigrid-based solvers.
For complex geometries, the use of finite element methods can also be advantageous~\cite{Jovanovic_2021}.

\begin{figure}
  \centering
  \includegraphics[width=7cm]{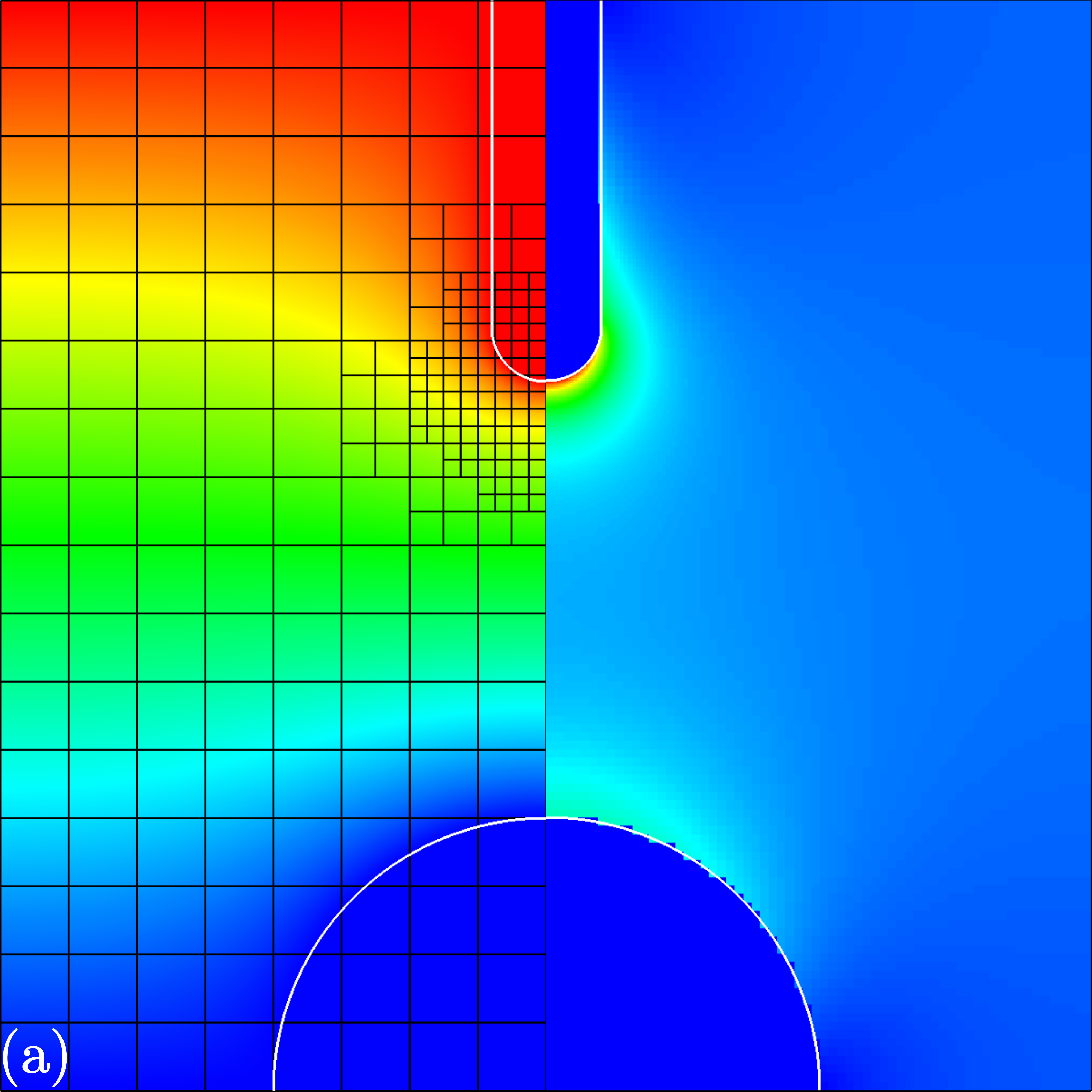}
  \includegraphics[width=8cm]{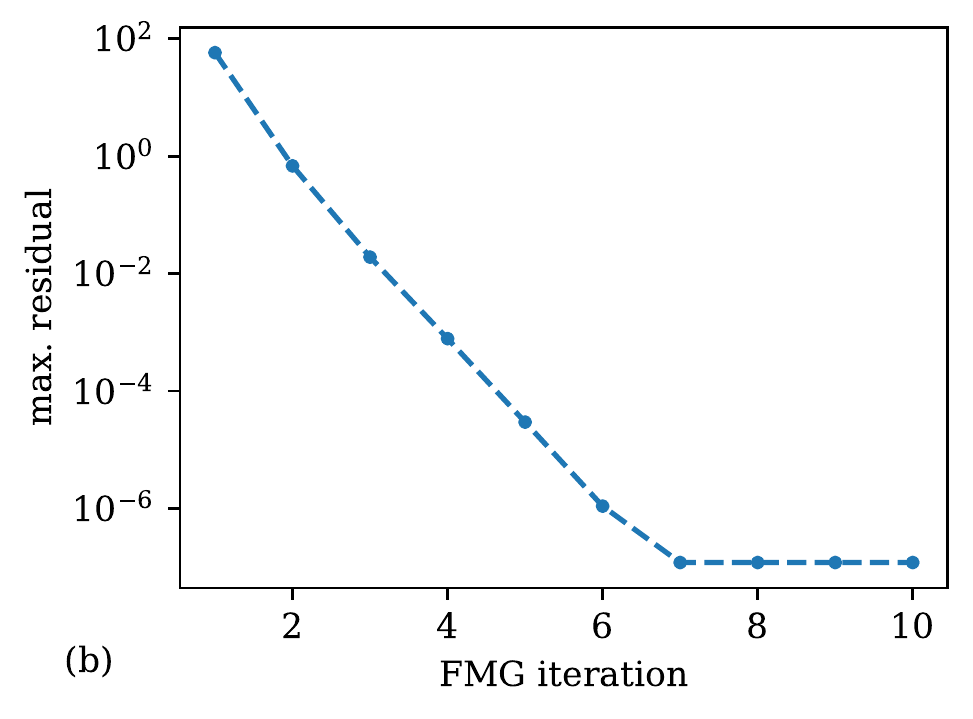}
  \caption{ a) The electric potential $\phi$ (left half) and the magnitude of its gradient $||\nabla \phi||$ (right half) for a two-electrode problem.
    The figure is a cross section of the 3D domain of size unity.
    Electrode contours are white, and the numerical mesh is illustrated, with each black square corresponding to a grid block of $8^3$ cells.
    b) Maximal residual versus FMG iteration.
  }
  \label{fig:two-electrodes}
\end{figure}

We consider a 3D geometry in which two electrodes are present, illustrated in figure~\ref{fig:two-electrodes}.
The computational domain is of size unity, and the following boundary conditions are used on its sides: $\phi = 1$ at the top, $\phi = 0$ at the bottom, and Neumann zero boundary conditions on the other sides.
A rod electrode is placed at the top, with radius $r_\mathrm{rod} = 0.05$, at which $\phi = 1$.
The corresponding level set function (in the top half) is obtained by computing the distance from a line segment, and then subtracting $r_\mathrm{rod}$.
On the bottom of the domain, a semi-sphere is placed with radius $r_\mathrm{sphere} = 0.25$, at which $\phi = 0$.

The numerical mesh has a spacing of $\Delta x = 1/512$ near the tip of the top electrode, and a resolution $\Delta x = 1/128$ elsewhere, as illustrated in figure~\ref{fig:two-electrodes}, which also shows the solution $\phi$, $||\nabla \phi||$ and the maximal residual versus FMG iteration.
The residual reduction factor is about 30-40 per FMG iteration.
The components of $\nabla \phi$ were computed on a staggered grid (on cell faces), taking the stored distances to boundaries into account.
For example, if there is a boundary between $\phi_{i,j,k}$ and $\phi_{i+1,j,k}$, then $\partial_x \phi$ at $i+1/2$ is approximated by $(\phi_b - \phi_{i,j,k})/(d \Delta x)$, where $\phi_b$ is the boundary value and $d$ is the relative distance to the boundary.
For cells whose center lies near the boundary but inside the electrodes, $||\nabla \phi||$ was set to zero.

\subsection{Computational cost}
\label{sec:computational-cost}

When new refinement is added to the mesh, the boundary detection method described in section~\ref{sec:boundary-detection} is performed.
This requires the evaluation of the numerical gradient of the LSF at every newly added grid cell\footnote{Note that in the majority of cases, the absence of an irregular boundary can be deduced from the parent grid, but some sharp features might only be detected on the fine grid.}.
Afterwards, the distance computation described in section~\ref{sec:lsf-distance} is performed for grid cells that are close to the boundary.
This requires a few tens of evaluations of the level set function per grid cell.
To keep these costs low, the LSF should be cheap to compute.
When removing refinement, no extra work is required.

In many applications, solutions have to be computed multiple times on the same numerical mesh, but with different right-hand sides.
The cost per multigrid iteration is then most important.
To illustrate these costs, we solve the test case with the spherical boundary described in section~\ref{sec:convergence-sphere} in 3D on a $512^3$ uniformly refined grid.
We consider block sizes of $8^3$, $16^3$ and $32^3$.
A smaller block size increases the adaptivity of the mesh, and it will reduce the total volume of grid blocks that intersect the boundary.
On the other hand, a smaller block size means that more blocks are required, leading to extra communication costs.

Table~\ref{tab:computational-cost} gives the time per FMG cycle in seconds for the various cases.
Note that the parallel scaling is not ideal.
The reason for this is that computations in a geometric multigrid method are relatively cheap, so that the speed with which data can be accessed from and written to memory is often the limiting factor.

\begin{table}
  \centering
  \begin{tabular}{l|lll}
     & $512^3/8^3$ & $512^3/16^3$ & $512^3/32^3$\\
    \hline
    4 cores & 5.55 & 2.97 & 2.36\\
    8 cores & 2.77 (100\%) & 1.53 (97\%) & 1.31 (90\%)\\
    16 cores & 1.62 (86\%) & 1.01 (74\%) & 0.87 (68\%)\\
    32 cores & 1.46 (48\%) & 0.95 (39\%) & 0.84 (35\%)
  \end{tabular}
  \caption{\jt{Computational time (in seconds) per FMG cycle for the 3D test case with a spherical boundary described in section~\ref{sec:convergence-sphere}, on a uniform grid with $512^3 \approx 134 \times 10^6$ cells.
    Results with block sizes of $8^3$, $16^3$ and $32^3$ are included, and parallel efficiencies compared to the case with 4 cores are indicated.
    The tests were performed using up to 32 cores of an AMD Epyc 7H12 processor.
    The computational times are an average over 40 consecutive FMG cycles.}}
  \label{tab:computational-cost}
\end{table}

\section{Conclusions}
\label{sec:conclusions}

We have presented a method to include irregular domain boundaries in a geometric multigrid solver.
The method was developed for quadtree/octree grids with adaptive refinement, using a cell-centered discretization, and it supports Dirichlet-type boundary conditions.
The location of boundary intersections is automatically determined from a level set function, which is to be provided as input.
For grid blocks near the interface, custom operator stencils are stored.
However, the computational cost of handling these custom blocks is comparable to that of regular blocks away from boundaries, and in both cases, point-wise multigrid smoothers are employed.
We have illustrated the numerical convergence of the method by considering spherical boundaries in the unit square and unit cube.
Furthermore, the robustness and computational efficiency of the method were examined with several test cases with sharp boundaries.\\

\bibliographystyle{elsarticle-num}
\bibliography{references}







\end{document}